\documentclass{article}
\usepackage{epsfig}
\usepackage{color}
\usepackage{amsmath,amssymb,amscd,amsthm} 
\usepackage{mathdots} 
 
\title{Analogues of Kahan's method for higher order equations of 
higher degree}  
\author{A.N.W. Hone\thanks{Work begun on leave at 
School of Mathematics \&  
Statistics, University of New South Wales, 
NSW 2052, Australia.}~\\
School of Mathematics, 
Statistics \& Actuarial Science~\\ 
University of
Kent~\\ Canterbury CT2 7NF, UK. \\ \\
 G.R.W. Quispel~\\ 
Department of Mathematics and Statistics~\\
La Trobe University~\\
Bundoora~\\ Victoria 3086,
Australia.
}

\newcommand{\br}{\begin{rem}}
\newcommand{\er}{\end{rem}}
\newcommand{\bex}{\begin{exa}}
\newcommand{\eex}{\end{exa}}
\newcommand{\bd}{\begin{Def}}
\newcommand{\ed}{\end{Def}}
\newcommand{\bt}{\begin{theorem}}
\newcommand{\et}{\end{theorem}}
\newcommand{\bl}{\begin{lemma}}
\newcommand{\el}{\end{lemma}}
\newcommand{\be}{\begin{equation}}
\newcommand{\ee}{\end{equation}}
\newcommand{\bea}{\begin{eqnarray}}
\newcommand{\eea}{\end{eqnarray}}

\newcommand{\adots}{\mathinner{\mkern2mu\raise1pt\hbox{.}\mkern2mu
\raise4pt\hbox{.}\mkern2mu\raise7pt\hbox{.}\mkern1mu}}

\newcommand{\beq}{\begin{equation}}  
\newcommand{\eeq}{\end{equation}}  
\newcommand{\bear}{\begin{array}}  
\newcommand{\eear}{\end{array}} 
\newcommand\la{{\lambda}}

\newcommand\al{{\alpha}}   
\newcommand\bet{{\beta}}   
\newcommand\om{{\omega}}

\newcommand\gam{{\gamma}}

\newcommand\si{{\sigma}}

\newcommand\rd{{\mathrm{d}}}  
\newcommand\ri{{\mathrm{i}}}  
\newcommand\tx{{\overline{x}}}
\newcommand\dx{{\underline{x}}} 
\newcommand\ty{{\overline{y}}}
\newcommand\bx{{\mathbf{x}}} 
\newcommand\by{{\mathbf{y}}} 
\newcommand\bb{{\mathbf{b}}} 
\newcommand\bof{{\mathbf{f}}}
\newtheorem{thm}{Theorem}[section] 

\newtheorem{propn}[thm]{Proposition} 

\newtheorem{rem}[thm]{Remark}

\newtheorem{exa}[thm]{Example}

\newenvironment{prf}{\trivlist \item [\hskip 
\labelsep {\bf Proof:}]\ignorespaces}{\qed \endtrivlist} 
 
\theoremstyle{remark}

\begin{document} 

\maketitle
 
\begin{abstract} 
Kahan introduced an explicit method of discretization for 
systems of first order differential equations with nonlinearities 
of degree at most two  
(quadratic vector fields).  
Kahan's method has attracted much interest due to the  fact 
that it preserves many of the geometrical properties of the 
original continuous system.  
In particular,  a large number 
of Hamiltonian systems of quadratic vector fields are known 
for which their Kahan discretization is a discrete integrable system. 
In this note, we introduce a special class of explicit 
order-preserving discretization 
schemes that are appropriate for certain systems of ordinary differential equations 
of higher order and higher degree.   \end{abstract}

\section{Introduction} 
\setcounter{equation}{0}

\noindent

Kahan's method is a special discretization scheme that provides 
an explicit method for integrating quadratic vector fields, given 
by systems of first order ordinary differential equations (ODEs) of the 
form 
\beq\label{1stodes}
\frac{\rd x_i}{\rd t} = f_i (x_1,\ldots,x_N), \qquad i=1,\ldots, N,  
\eeq  
where each function $f_i$ is a polynomial of total degree two in 
the independent variables $x_1,\ldots,x_N$ 
(see \cite{kahanlect, kahan1}). In order to specify 
Kahan's method, one should replace each derivative on the 
left-hand side of (\ref{1stodes}) 
by the forward difference, so 
that 
$$ 
\frac{\rd x_i}{\rd t}\rightarrow \Delta x_i:=\frac{\tx_i-x_i}{h}, 
$$ 
while terms of degrees two, one and zero appearing in each $f_i$ 
on the right-hand side are 
replaced according to the 
rules 
\beq\label{rules}
x_jx_k \rightarrow \frac{1}{2}(\tx_jx_k + x_j\tx_k), 
\quad x_j\rightarrow \frac{1}{2}(x_j+\tx_j), \quad 
c\rightarrow c , 
\eeq 
where $h$ is the time step and ${\tx}_i$ denotes the approximation 
to $x_i(t+ h)$. 

It was noticed some time ago that 
Kahan's method provides an effective integration scheme 
for the classic two-species Lotka-Volterra model 
\beq\label{lv} 
\frac{\rd x}{\rd t}  =  \alpha x(1-y), \qquad
\frac{\rd y}{\rd t} = y(x-1) 
\eeq
(with $\al>0$ being an arbitrary parameter), retaining the qualitative features of the orbits 
of the continuous system, namely the stability of orbits 
around the elliptic fixed point at $(x,y)=(1,1)$.
This was subsequently explained by the fact that 
the Kahan discretization of  (\ref{lv}), given by 
$$ \begin{array}{rcl}
(\tx - x)/h & = & \frac{\al}{2}\Big( x(1-\ty)+\tx(1-y)\Big), \\ 
 (\ty - y)/h & = & \frac{1}{2}\Big( y(\tx-1)+\ty(x-1)\Big),
\end{array}  
$$ 
preserves the same symplectic form 
$$ 
\om = \frac{\rd x\wedge \rd y}{xy}
$$ 
as the original Hamiltonian system \cite{sanzserna}. 
In the context of Lotka-Volterra models, a variant of 
Kahan's method with similar properties 
was discovered  by Mickens \cite{mickens},  who had previously considered 
various examples of nonstandard discretization methods \cite{mickensbook}, 
but a more  rapid growth of interest in Kahan's method 
began when Hirota and Kimura independently proposed      
the rules (\ref{rules}) for the 
discretization of the Euler equations for rigid body motion, 
finding that the resulting discrete system is 
also completely integrable \cite{hirotakimura}, and this has 
led to the search for other discrete integrable systems arising in this way 
\cite{hp, pps, pz, pss}.

Many of the geometrical properties of Kahan's method for 
quadratic vector fields are based 
on the polarization identity for quadratic forms \cite{cmoq},  
and recently this has led to 
a generalization of Kahan's method that can cope with vector fields of 
degree three or more, by using higher degree analogues of polarization
\cite{cmmoq}.  One disadvantage of the latter method 
for higher degree vector fields is that, 
in common with multistep  methods in numerical analysis, 
one must use extra grid points for the discretization, so the original ODE system 
does not provide enough initial values to start the iteration of the discrete version. 
However, 
if one is looking for a discretization scheme that preserves integrability 
or other geometric properties of  
ODEs, then it is desirable for  the initial value space 
of the discrete system to have 
the same dimension as that of the continuous one.        
Here we would  like to suggest a discretization scheme 
with the latter property, which is a natural generalization of Kahan's method to higher order and higher 
degree.

The idea is to consider a system of ODEs of order $n\geq 1$, with the right-hand 
sides being functions of the coordinates $x_1,\ldots, x_N$ only, of the form 
\beq\label{odes}
\frac{\rd^n x_i}{\rd t^n} = f_i (x_1,\ldots,x_N), \qquad i=1,\ldots, N,  
\eeq  
where each function $f_i$ is a polynomial of maximal degree $n+1$. For $n=1$ this 
is a quadratic vector field, which one can discretize using Kahan's method. 
In the next section, we present an explicit discretization scheme for 
systems of the form (\ref{odes}), valid for any $n\geq1$, 
which reduces to Kahan's method when $n=1$. The first new case is $n=2$, 
corresponding to systems of Newton equations, which are relevant in many
applications. 
We illustrate this in section 3 by considering the discretization of the 
motion of a single particle moving in a quartic potential. The latter is 
one of the simplest examples of an integrable Hamiltonian system, and 
it turns out that the discrete version produced by the method is also integrable, 
with a conserved quantity and an invariant symplectic form. In section 4, 
we consider a different  example of fourth order, namely a nonlinear beam equation, 
and briefly compare the discretization obtained by the new method with 
another discretization obtained by applying an approach similar to Kahan's 
directly to the Lagrangian of the continuous system.

\section{A higher order version of Kahan's method}

\setcounter{equation}{0}

For $n=2$, 
(\ref{odes}) becomes a system of Newton equations, assumed to have  polynomial forces of degree at most three, which 
can be conveniently written as 
\beq\label{newton} 
\frac{\rd^2 x_i}{\rd t^2} = \sum_{0\leq j_1\leq j_2\leq j_3\leq N} c_{i,j_1,j_2,j_3} \, x_{j_1}x_{j_2}x_{j_3}, \qquad i=1,\ldots, N,
\eeq 
where $c_{ijk\ell}$ are arbitrary coefficients, and 
we have included an additional dummy variable $x_0=1$ to allow terms of degree less than three 
to be included within the same summation. Then 
to discretize (\ref{newton}) we propose the 
following: 
\beq\label{cubic}
\frac{{\tx}_i - 2x_i + {\dx}_i}{h^2} = \frac{1}{6}\sum_{\si\in S_3}\sum_{0\leq j_1\leq j_2\leq j_3\leq N}
c_{i,j_1,j_2,j_3} \, {\dx}_{j_{\si (1)}}x_{j_{\si (2)}}{\tx}_{j_{\si (3)}}, \qquad i=1,\ldots, N;  
\eeq 
the first summation is over  permutations $\si$ in $S_3$, the symmetric group on three symbols, 
and ${\tx}_i=x^{(1)}$, ${\dx}_i=x^{(-1)}$ are the approximations to $x_i(t\pm h)$, with time step $h$.   
For terms of degree three, with each variable $x_j$ appearing at the three adjacent lattice points
$\dx_j=x^{(-1)}_j$, $x_j=x^{(0)}_j$, $\tx_j =x^{(1)}_j$,  
  the replacement rule is described explicitly by 
\beq\label{cube}
x_jx_kx_\ell \rightarrow \frac{1}{6}\Big(\dx_j x_k \tx_\ell+\dx_j \tx_k x_\ell 
+x_j\dx_k\tx_\ell + x_j \tx_k\dx_\ell 
+\tx_j x_k\dx_\ell +\tx_j\dx_kx_\ell\Big), 
\eeq
while for terms of degree two the rule is obtained by setting $\ell = 0$, so that $x_\ell\to x_0=1$ in the above,
and for the linear terms one can set $k=\ell=0$, so that the rule for terms of degree less than three is 
\beq\label{quads}
x_jx_k \rightarrow \frac{1}{6}\Big(
\dx_j x_k +\dx_j \tx_k  
+x_j\dx_k + x_j \tx_k 
+\tx_j x_k +\tx_j\dx_k
\Big), \quad 
x_j\rightarrow 
 \frac{1}{3}\Big(
\dx_j    
+x_j
+\tx_j 
\Big), \quad c\rightarrow c.
\eeq 
  
Following the approach of \cite{cmmoq}, a second order system of equations 
can be 
written in vector form as  
\beq\label{2ndorder}
\frac{\rd^2 \bx}{\rd t^2}=\bof (\bx),
\eeq 
where each component of the vector of functions $\bof = (f_1,f_2,\ldots,f_N)^T$ is a polynomial of degree at most three, 
and then the replacement rules (\ref{cube}) and (\ref{quads}) are equivalent to the formula 
\beq\label{vecf}
\frac{\overline{\bx} - 2\bx +\underline{\bx}}{h^2} = 
\frac{9}{2}
\bof\Big(
\frac{\overline{\bx} +\bx+\underline{\bx}}{3}
\Big)
-\frac{4}{3}
\left( 
\bof\Big( 
\frac{\underline{\bx} + \bx}{2}
\Big)+\bof\Big( 
\frac{\underline{\bx} + \overline{\bx}}{2}
\Big)+\bof\Big( 
\frac{\bx+\overline{\bx}}{2}\Big)
\right)
+\frac{1}{6}\Big( \bof( \underline{\bx} ) +\bof( {\bx} )+ \bof( \overline{\bx} )
\Big).
\eeq 

\begin{propn} 
The discretization (\ref{vecf})
commutes with affine transformations 
\beq\label{affine} 
\by \mapsto \bx =A \by + \bb, 
\eeq 
where $A\in GL(N,\mathbb{R}) $ is a 
constant matrix and $\bb\in\mathbb{R}^N$ is a 
vector of constants. 
\end{propn} 
\begin{prf} Under the transformation (\ref{affine}),
$\bof (\cdot) $ in (\ref{2ndorder})
is replaced by $A^{-1} \bof (A\cdot +\bb )$. 
Upon substituting (\ref{affine}) 
and its shifted versions into (\ref{vecf}),
it is not hard to check that the same occurs for each appearance of $\bof$ 
on the right-hand side. \end{prf} 

The symmetric  replacement rules above generalize to any order $n\geq 1$, so that 
for a system of $n$th order ODEs (\ref{odes}) with  right-hand sides all of degree $n+1$  the discretization  becomes
\beq\label{symn} 
\Delta^n x_i=
\frac{1}{(n+1)!}
\sum_{\si\in S_{n+1}}\sum_{0\leq j_1\leq \cdots \leq j_{n+1}\leq N}
c_{i,j_1,\ldots,j_{n+1}} \, x_{j_{\si (1)}}x_{j_{\si (2)}}^{(1)}x_{j_{\si (3)}}^{(2)}\cdots x_{j_{\si (n+1)}}^{(n)},  
\quad i=1,\ldots, N, 
\eeq 
with $x_i^{(1)}={\tx}_i, x_i^{(2)}=\overline{\tx}_i, \ldots , x_j^{(n)}$ corresponding to shifts by steps of $h, 2h, \ldots, nh$. 
On the left-hand side of (\ref{symn}) we have replaced the $n$th derivative by the $n$th power of the forward difference operator, and for convenience 
we have written everything on the right-hand side in terms of forward shifts of the variables $x_j$. 
The discretization (\ref{symn}) reduces to Kahan's method when $n=1$, and to (\ref{cubic}) when $n=2$,  
modulo shifting the lattice points $-1,0,1$ in the latter up to $0,1,2$.

Clearly there are other 
choices of discrete  $n$th derivative that one could take, and other affine combinations of terms with the same homogeneous degree could be chosen 
while preserving the continuum limit. We have taken the most symmetrical choice in (\ref{symn}), because it is manifestly linear in each of 
the highest shifts $x_1^{(n)},x_2^{(n)}\ldots,x_N^{(n)}$, so it can be explicitly solved for each of these quantities to yield 
rational functions of all the lower shifts. It is also linear in each of the lowest shifts $x_i=x_i^{(0)} $ for $i=1,\ldots,N$, so 
it can be explicitly solved for these as well. Thus 
(\ref{symn}) is an implicit way of writing an explicit birational map in dimension $nN$, corresponding to 
\beq\label{birat}
\begin{array}{rl} 
(x_1^{(0)},\ldots,x_N^{(0)},\,\, x_1^{(1)},\ldots,x_N^{(1)},\,\,\,\ldots & , x_1^{(n-1)},\ldots,  x_N^{(n-1)}) 
\\ 
& \mapsto 
(x_1^{(1)},\ldots,x_N^{(1)},\,\,x_1^{(2)},\ldots,x_N^{(2)},\,\,\,\ldots\,\,\,, x_1^{(n)},\ldots,  x_N^{(n)}) .
\end{array}
\eeq

\section{Discretization of a quartic oscillator}
\setcounter{equation}{0}

To see why it might be worth investigating these higher Kahan-like schemes, 
we start by presenting the following 
 example: $n=2$ with a cubic force on a particle in one dimension, generated by a 
natural Hamiltonian with a quartic potential, that is 
$$ 
H=\frac{1}{2}p^2 +\frac{1}{4}ax^4 +\frac{1}{3}bx^3+\frac{1}{2}cx^2+dx,
$$   
which yields the Newton equation 
\beq\label{newt} 
\ddot{x} = -ax^3-bx^2-cx-d.
\eeq
This is an integrable system par excellence, and the generic  level sets $H=\,$const are quartic curves of genus one
in the $(x,p)$ plane.   
The discretization (\ref{cubic}) applied to (\ref{newt}) produces a difference equation of second order, given by 
\beq\label{qrt} 
{\tx} = \frac{(3-\gam) x -\delta- (\bet x +\gam) {\dx}} {\bet x +\gam+(\al x+\bet ){\dx}} , 
\eeq 
where 
$$
\al=ah^2, \quad \bet =\frac{bh^2}{3}, \quad \gam =1+\frac{ch^2}{3}, \quad \delta =dh^2.
$$

The map (\ref{qrt}) is an  example of a QRT map \cite{QRT1}, but let us suppose that 
we do  not know the geometric properties of this map. To find these properties, 
such as the existence of a preserved measure, and first and second 
integrals of the map (\ref{qrt}), we will look for preserved Darboux 
polynomials, as detailed in our recent work \cite{cemoqtk, cemoqt}. 
To this end, 
we write the second order equation (\ref{qrt}) as two first order ones, namely 
\beq\label{1stord}
\tx = y, 
\qquad
\ty = \frac{(3-\gam) y-\delta - (\bet y +\gam)x} {\bet y+\gam  +(\al y+\bet )x }  ,  
\eeq 
and look for polynomials $P$ satisfying
\beq\label{darboux} 
P(\tx, \ty )=J(x,y) P(x,y), 
\eeq 
where $J$ is the Jacobian determinant of the map (\ref{1stord}), i.e.  
\beq \label{jac}
J(x,y )=\frac{(\bet y +\gam)^2+(\al y+\bet)\Big((3-\gam)y-\delta\Big) }{(\al xy +\bet (x+y)+\gam)^2} .
\eeq 
Substituting (\ref{jac}) into (\ref{darboux}), and looking for polynomials up to total degree four 
in $x$ and $y$, we find two linearly independent solutions, given by 
$$  
P_1=\al xy +\bet  (x+y)+\gam, 
\quad
P_2= (\al \gam-\bet^2) x^2y^2+ \epsilon xy(x+y) +\zeta (x^2+y^2)-(3-\gam)^2 xy
+(3-\gam) \delta(x+y)-\delta^2, 
$$  
with 
$$ 
\epsilon = \al \delta +\bet (3-\gam), \qquad 
\zeta = 
\bet\delta+\gam(3-\gam).
$$
It follows that the map (\ref{qrt}) is measure-preserving, with the invariant symplectic form 
\beq\label{meas} 
\frac{\rd x\wedge \rd y}{P_1}= \frac{\rd x \wedge \rd y}{\al xy + \bet (x+y)+\gam}, 
\eeq 
and the first integral 
\beq\label{firstint}
\frac{P_2}{P_1}= \frac{
(\al \gam-\bet^2) x^2y^2+ \epsilon xy(x+y) +\zeta (x^2+y^2)-(3-\gam)^2 xy
+(3-\gam) \delta(x+y)-\delta^2}
{\al xy +\bet  (x+y)+\gam}.
\eeq 
Hence the integrability is preserved by the discretization in this case, and we recover the standard property of a QRT map, that 
it preserves a pencil of biquadratic curves, here given by 
$$ 
 \lambda P_1(x,y) +  P_2(x,y)=0.  
$$
Moreover, in the continuum limit $h\to 0$, the standard area form $\rd x\wedge \rd y$ and the Hamiltonian $H$ are recovered from (\ref{meas}) and (\ref{firstint}) respectively, since from  
$
y=x+hp+O(h^2)
$ 
we find 
$$ 
P_1 = 1+O(h^2), \qquad P_2= 4Hh^2+O(h^3).
$$ 

The equation (\ref{newt}) includes Duffing's equation, which is the case $b=d=0$, and also 
the second order ODE for the Weierstrass $\wp$ function, which arises when $a=c=0$. 
 In \cite{potts1}, another replacement rule is used for the cubic and linear terms in Duffing's equation, 
somewhat less symmetrical than the one defined by (\ref{cube}), and it is shown that if the 
coefficients and denominator in the second difference operator  are replaced by suitable functions of the parameters 
and the time step $h$ then this alternative rule results in a discretization that is exact, in  the sense that the iterates of the difference equation 
interpolate the solution of the original ODE. Similarly, in \cite{potts2} an exact discretization is obtained for the case 
corresponding to the Weierstrass $\wp$ function, with only quadratic and constant terms on the 
right-hand side. However, in the latter case,  the exact discretization (derived from the addition 
formula for the $\wp$ function)  requires not only a different 
replacement rule for the quadratic terms compared with (\ref{quads}), but also extra cubic and linear 
terms that must be included, with a coefficient which is $O(h^2)$. 
When $a=c=0$, the equation (\ref{newt}) can be rewritten as a quadratic 
vector field, namely 
$$ 
\frac{\rd x}{\rd t}=p, \qquad 
\frac{\rd p}{\rd t}= -bx^2 -d, 
$$ 
so that Kahan's method can be applied, as in \cite{pps}, resulting in 
a first order discrete system which is equivalent to  a 
second order difference equation for $x$, namely 
\beq\label{qrtadd}
\tx +\dx = \frac{4x-2\delta}{3\beta x+2} 
\eeq (where we set $\beta =bh^2/3$, $\delta=dh^2$ as before). 
The equation (\ref{qrtadd}) is a QRT map in additive form, clearly 
of a different type to (\ref{qrt}), which becomes 
$$ 
{\tx} = \frac{2 x -\delta- (\bet x+1)  {\dx}} {\bet (x +{\dx})+1}
$$ 
when $\al =0$,  $\gam =1$. To see that they are really different 
QRT maps, in the sense that they are not related to one another via  so-called curve-dependent McMillan maps \cite{ir}, 
observe that the pencil of invariant biquadraic curves corresponding to (\ref{qrtadd}) is 
$$\lambda -\beta^2x^2y^2+\frac{4}{3}\beta xy(x+y)+\frac{4}{3}(x^2+y^2)-\frac{2}{3}(4+\beta \delta)xy+\frac{4}{3}\delta (x+y)=0,$$
whereas when $\alpha=0$, $\gamma=1$ the pencil   $\lambda P_1(x,y) +  P_2(x,y)=0$  for (\ref{qrt}) reduces to  one of a different type, namely
$$
\lambda\Big(1+\beta (x+y) \Big)- \beta^2x^2y^2+2\beta xy(x+y)+(\beta \delta+2)(x^2+y^2)-4xy+2\delta (x+y)-\delta^2
=0. 
$$

\section{Two discretizations of a nonlinear beam equation} 
\setcounter{equation}{0}

Vibrating beams were considered by Leonardo da Vinci \cite{davinci}, but the 
traditional 
theory of vibrations of a beam is usually attributed to Euler and Bernoulli 
\cite{beam}, being described by a partial differential equation (PDE) of fourth order, 
which in dimensionless form is given by 
$$ 
\frac{\partial^2w}{\partial t^2} +
\frac{\partial^4w}{\partial x^4} =Q
.
$$
 For the case of a static beam, the equation has the form 
\beq 
\label{beameq}
\frac{\rd^4 w}{\rd x^4}=Q, 
\eeq
where $w=w(x)$ is the vertical deflection of the beam, which lies horizontally 
along the $x$-axis. The standard beam model is linear, with the 
distributed load $Q$ on the 
right-hand size being a constant (or more generally, a function of $x$, 
the independent variable). 
However, here we consider a more general nonlinear version of the model,  
derived from a second order Lagrangian of the form 
\beq\label{beamlag}
L = \frac{1}{2}\left(\frac{\rd^2 w}{\rd x^2}\right)^2 -V(w), 
\eeq 
which gives a nonlinear load function 
$$ 
Q(w)=\frac{\rd V}{\rd w}.
$$ 
In the linear case, the model was considered recently from 
the viewpoint of a Hamilton-Jacobi approach to higher order 
implicit systems \cite{sardon}, while a coupled PDE system of beam equations with cubic 
nonlinearity was analysed in \cite{shixu}. From the second order Lagrangian 
(\ref{beameq}), 
we can introduce the Ostrogradsky variables (see \cite{blaszak}, 
for instance), given by   
$$ 
q_1 = w, \quad q_2 = w' , 
\quad
p_1 =\frac{\partial L}{\partial w'} -\frac{\rd}{\rd x} \left( 
\frac{\partial L}{\partial w''}\right)=-w''', 
\quad 
p_2 = \frac{\partial L}{\partial w''}  =w'', 
$$ 
where the primes denote derivatives with respect to the 
independent variable $x$. 
Then $(q_1,p_1)$, $(q_2,p_2)$ provide two pairs of 
canonically conjugate positions and momenta, and the Euler-Lagrange 
equation 
\beq\label{el} 
\frac{\rd^2}{\rd x^2} \left( 
\frac{\partial L}{\partial w''}\right) - 
\frac{\rd}{\rd x} \left( 
\frac{\partial L}{\partial w'}\right) 
+ 
\frac{\partial L}{\partial w}=0, 
\eeq  
which for the Lagrangian (\ref{beamlag}) is given by (\ref{beameq}) 
with $Q=\rd V/\rd w$, is equivalent to 
Hamilton's equations for the 
Hamiltonian function 
$$ 
H = \frac{1}{2}(p_2)^2 +q_2 p_1 +V(q_1). 
$$ 
For the sake of concreteness, we consider the case of an odd 
potential 
$$ 
V(w) = \frac{a}{5} w^5 + \frac{b}{3} w^3 + c w, 
$$ 
so that the nonlinear beam equation is given by 
\beq\label{nlbeam} 
w'''' = aw^4 + bw^2+c. 
\eeq 

To begin with, we consider the result of applying the discretization 
rule (\ref{symn}) to the nonlinear beam equation (\ref{nlbeam}), 
which produces a difference equation of fourth order, of the form 
\beq\label{dbeam} 
\Delta^4 w = F(w^{(0)},w^{(1)}, w^{(2)},w^{(3)},w^{(4)})
\eeq 
for a function $F$ that is a sum of terms of total degree four, two and zero.
This can be written more symmetrically by shifting down by two steps, 
to yield 
\beq\label{dbeam_1} 
\frac{w^{(-2)}-4w^{(-1)}+6w^{(0)}-4w^{(1)}+w^{(2)}}{h^4} 
= F_{4}+F_{2} +c, 
\eeq 
where the quartic terms are 
\beq\label{quart} 
F_4=\frac{a}{5}\left(\begin{array}{c}w^{(-2)}w^{(-1)}w^{(0)}w^{(1)} 
+w^{(-2)}w^{(-1)}w^{(0)}w^{(2)}+
w^{(-2)}w^{(-1)}w^{(1)}w^{(2)} \\
+
w^{(-2)}w^{(0)}w^{(1)}w^{(2)}+ 
w^{(-1)}w^{(0)}w^{(1)}w^{(2)}\end{array}\right),  
\eeq 
and the quadratic terms are given by 
\beq\label{quadr} 
F_2=\frac{b}{10}\left(\begin{array}{c}w^{(-2)}w^{(-1)}+w^{(-2)}w^{(0)} 
+w^{(-2)}w^{(1)}+w^{(-2)}w^{(2)}+
w^{(-1)}w^{(0)}   \\ 
+w^{(-1)}w^{(1)}+ 
w^{(-1)}w^{(2)}+w^{(0)}w^{(1)}+ 
w^{(0)}w^{(2)}+w^{(1)}w^{(2)}\end{array}\right).   
\eeq 

It turns out that the birational  map defined by (\ref{dbeam_1}) is measure-preserving. This is a consequence of the fact that 
the formula for the right-hand side of  (\ref{dbeam}) is both linear and symmetric in its arguments, so that the derivatives with respect to the 
highest and lowest shifts, 
namely
\beq \label{fder}
\frac{\partial F}{\partial w^{(0)}}=G(w^{(1)},w^{(2)},w^{(3)},w^{(4)}) ,
\qquad 
\frac{\partial F}{\partial w^{(4)}}=H(w^{(0)},w^{(1)},w^{(2)},w^{(3)}) ,
\eeq
are very closely related to one another.
\begin{propn}
The discretization (\ref{dbeam_1}) preserves the volume form 
$$
\Omega = \frac{1}{1-h^2H(w^{(-2)},w^{(-1)},w^{(0)},w^{(1)}) }\, \rd w^{(-2)}\wedge \rd w^{(-1)}\wedge \rd w^{(0)}\wedge \rd w^{(1)},
$$
where $H$ is defined by (\ref{fder}).
\end{propn}
\begin{prf}
Upon taking the differential of both sides of (\ref{dbeam_1}) , we obtain the equation 
$$
\Big(1-h^2 G(w^{(-1)},w^{(0)},w^{(1)},w^{(2)}) \Big)\rd w^{(-2)} + 
\Big(1-h^2 H(w^{(-2)},w^{(-1)},w^{(0)},w^{(1)}) \Big)\rd w^{(2)} + \cdots=0,
$$
where the ellipsis denotes terms that are linear in $\rd w^{(-1)}$, $ \rd w^{(0)}$ and $ \rd w^{(1)} $. The result then follows
from taking the wedge product of the above with $\rd w^{(-1)}\wedge \rd w^{(0)}\wedge \rd w^{(1)} $, 
and noting the identity 
$$
 G(w^{(-1)},w^{(0)},w^{(1)},w^{(2)})=H(w^{(-1)},w^{(0)},w^{(1)},w^{(2)}), 
$$
which follows from the symmetry of $F$.
\end{prf}

 When $a,b$ are not both zero,  so that the nonlinear terms are present, the above discretization of (\ref{nlbeam}) cannot be obtained from a second order discrete Lagrangian 
of the form 
$$
{\cal L}
={\cal L}(w^{(n)},
w^{(n+1)},
w^{(n+2)}),
$$
since the discrete Euler-Lagrange equations 
\beq\label{del}
\sum_{i=0}^2\frac{\partial}{\partial 
w^{(n)}}
{\cal L}(w^{(n-i)},
w^{(n+1-i)},
w^{(n+2-i)})=0
\eeq
do not generate terms containing products $w^{(j)} w^{k)}$ with $|j-k|>2$. In order to obtain a discretization 
with a Lagrangian structure, we fix $n=0$ and take a discrete Lagrangian of the form 
$${\cal L}(w^{(0)},
w^{(1)},
w^{(2)})={\cal T}-{\cal V}, 
$$
where 
the discrete fourth derivative is generated by the term 
$${\cal T}=\frac{1}{2h^4}\Big(2(w^{(0)}-w^{(1)})^2 - (w^{(0)}-w^{(2)})^2+
2(w^{(1)}-w^{(2)})^2\Big),  
$$
and the other terms are specified by 
$${\cal V}={\cal V}_5+{\cal V}_3+\frac{c}{3}\Big(w^{(0)}+
w^{(1)}+
w^{(2)}\Big)
$$
with 
\beq\label{v5} 
{\cal V}_5=\frac{a}{5}
\left(\begin{array}{c} 
\al_0
w^{(0)}(w^{(1)})^3 w^{(2)}+
\frac{1}{2}\al_1 \Big(
(w^{(0)})^2 (w^{(1)})^3
+
(w^{(1)})^2 (w^{(2)})^3 \Big) 
\\ 
+
\frac{1}{2}\al_2 \Big(
(w^{(1)})^2 (w^{(0)})^3
+
(w^{(2)})^2 (w^{(1)})^3 \Big) 
+
\frac{1}{2}\al_3 \Big(
w^{(0)} (w^{(1)})^4
+
w^{(1)} (w^{(2)})^4 \Big)
\\ 
+
\frac{1}{2}\al_4 \Big(
w^{(1)} (w^{(0)})^4
+
w^{(2)} (w^{(1)})^4 \Big)
+
\frac{1}{3}\al_5 \Big(
(w^{(0)})^5
+ 
(w^{(1)})^5
+
(w^{(2)})^5
\Big)
\end{array} 
\right),
\eeq 
and 
\beq\label{v3} 
{\cal V}_3=\frac{b}{3}\left(
\begin{array}{c}
\bet_0
w^{(0)}w^{(1)} w^{(2)}+
\frac{1}{2}\bet_1 \Big(
w^{(0)} (w^{(1)})^2
+
w^{(1)} (w^{(2)})^2 \Big) 
+
\frac{1}{2}\bet_2 \Big(
w^{(1)} (w^{(0)})^2
+
w^{(2)} (w^{(1)})^2 \Big)
\\ 
+
\frac{1}{3}\bet_3 \Big(
(w^{(0)})^3
+ 
(w^{(1)})^3
+
(w^{(2)})^3
\Big)\end{array}
\right), 
\eeq 
where, in (\ref{v5}) and (\ref{v3}) we 
have taken affine combinations, 
so that the coefficients are required to satisfy 
$$ 
\sum_{j=0}^5 \al_j=1 = \sum_{j=0}^3 \bet_j 
$$ 
in order 
to ensure the correct continuum limit, 
and we have included all 
possible terms of degrees 5 and 3, respectively, 
except those whose discrete variational derivative 
produces expressions of degree greater than one in 
$w^{(-2)}$ or $w^{(2)}$ (we have also grouped 
together terms having the same variational derivative). 
Hence we arrive at a discretization of (\ref{nlbeam}) 
which is explicit and birational, being given by 
\beq\label{dbeam_2} 
\frac{w^{(-2)}-4w^{(-1)}+6w^{(0)}-4w^{(1)}+w^{(2)}}{h^4} 
= \hat{F}_{4}+\hat{F}_{2} +c, 
\eeq 
where the quartic and  quadratic terms are given by  
\beq\label{lquart} 
\hat{F}_4=\frac{a}{5}\left(\begin{array}{l}
\al_0 \Big( w^{(-2)}(w^{(-1)})^3 
+3 w^{(-1)}(w^{(0)})^2 w^{(1)}
+(w^{(1)})^3 w^{(2)}\Big) \\
+  \al_1 \Big(
3 (w^{(-1)})^2(w^{(0)})^2
+2 w^{(0)}(w^{(1)})^3
\Big) 
+  \al_2 \Big(
2 (w^{(-1)})^3 w^{(0)}
+3 (w^{(0)})^2(w^{(1)})^2
\Big) \\ 
+  \al_3 \Big(
4w^{(-1)}(w^{(0)})^3 
+ (w^{(1)})^4
\Big) 
+  \al_4 \Big(
(w^{(-1)})^4
+4(w^{(0)})^3 w^{(1)}
\Big)
+5\al_5(w^{(0)})^4
\end{array}\right),  
\eeq 
\beq\label{lquadr} 
\hat{F}_2=\frac{b}{3}\left(
\begin{array}{c}
\bet_0\Big(w^{(-2)}w^{(-1)}+w^{(-1)}w^{(1)}+ w^{(1)}w^{(2)}\Big) 
 +
 \bet_1\Big(2w^{(-1)}w^{(0)}+(w^{(1)})^2\Big)
 \\ 
 +
 \bet_2\Big((w^{(-1)})^2+ 2w^{(0)}w^{(1)}\Big) 
 +
 3\bet_3(w^{(0)})^2
\end{array}
\right),   
\eeq 
respectively. 
A general approach to Lagrangian fourth-order difference equations and their continuum limits appears in the recent paper \cite{gubb}.

An advantage of using the Lagrangian discretization (\ref{dbeam_2}) 
is that it is symplectic; so it is a birational symplectic integrator. This can be seen from  
the discrete analogue of the Ostrogradsky transformation, 
introduced in \cite{bruschi}, which provides canonical variables 
$q_1,p_1,q_2,p_2$ via the formulae
\beq\label{canonical} 
q_1=w^{(0)}, \, q_2 = w^{(1)}, \,
p_1 = {\cal L}_1(w^{(-1)}, w^{(0)},w^{(1)}) + 
{\cal L}_2(w^{(-2)}, w^{(1)},w^{(0)}), 
\,
p_2 = {\cal L}_2(w^{(-1)}, w^{(0)},w^{(1)}), 
\eeq  
where 
$$ 
{\cal L}_j =\frac{\partial {\cal L}}{\partial w^{(j)} }
( w^{(0)},w^{(1)},w^{(2)}), \qquad j=0,1,2. 
$$ 
In terms of these variables, the four-dimensional map 
defined by (\ref{dbeam_2}) preserves the canonical 
symplectic form 
$$ 
\omega =
\rd p_1\wedge \rd q_1 + \rd p_2 \wedge \rd q_2, 
$$ 
and this immediately implies that it preserves the volume form
$\omega\wedge\omega$, so it is measure-preserving. 

Qualitatively it appears that the approximate solutions 
of (\ref{nlbeam}) provided by these two discretizations 
are somewhat similar. To see this, one can consider 
solutions in the neighbourhood of a fixed point. 
If  $ab\neq 0$ then,  by scaling $w$ and $x$, the 
parameters can be taken as 
$$ 
a=1, \qquad b=-2\epsilon, \qquad c=1-\delta, 
$$ 
with $\epsilon^2=1$ and $\delta$ arbitrary. 
Then (\ref{nlbeam}) has fixed 
points at $w=\pm \sqrt{\epsilon\pm \sqrt{\delta}}$, so that  
$\delta\geq 0$ is a necessary condition for reality, and then 
generically there are either four, two or zero real fixed 
points depending on the choice of $\epsilon=\pm 1$ and the value of 
$\delta$. 
In 
particular,   
let us take the case 
$$ 
\epsilon=1, \qquad 0<\delta<1 
$$ 
when there are four real fixed points, one of which is at $w=w^*$, 
where 
$$ w^* =  \sqrt{1+ \sqrt{\delta}}.$$
The eigenvalues of the linearization of   (\ref{nlbeam}) 
around this point consist of a real pair $\pm \gamma$ and an imaginary 
pair $\pm \ri \gamma$, for $\gam = (4w^*\sqrt{\delta})^{1/4}$, 
corresponding to one stable direction, one unstable direction, and
a two-dimensional centre manifold. The  discretizations (\ref{dbeam_1}) 
and (\ref{dbeam_2}) both have the same fixed points as the 
original differential equation, and using the fact that 
(\ref{dbeam_1}) is reversible, and that (\ref{dbeam_2}) is symplectic (and also 
reversible), together with standard facts about linear stability of 
reversible/symplectic maps \cite{hm, lbr}, in each case the 
characteristic polynomial 
of the linearization around a 
fixed point is palindromic (equivalently,  $\la$ is a root if 
and only if $\la^{-1}$ is). If we consider the linearization 
around $w^*$, then in both cases we find two real eigenvalues that are 
reciprocals of one another, corresponding to the stable and unstable directions, 
together with a complex conjugate pair of eigenvalues of modulus one, giving 
a two-dimensional centre manifold, just as for the differential equation;
and similar considerations apply to the other fixed points. Thus, to a 
first approximation, the qualitative behaviour 
of the two discretizations is the same.  

\section{Conclusions}

We have found that the higher order analogue of Kahan's method proposed here preserves integrability in the  second order example of the quartic oscillator (\ref{newt}) that we have considered, 
while in the case of a nonlinear beam equation of fourth order the resulting discretization (\ref{dbeam_1}) is measure-preserving, and its qualitative behaviour 
looks similar to that of the Lagrangian discretization (\ref{dbeam_2}). In future work we would like to apply this discretization method to other ODE systems of higher order, as well as 
looking for first integrals of the particular fourth order maps (\ref{dbeam_1}) and (\ref{dbeam_2}) using the method of discrete Darboux polynomials as in \cite{cemoqtk, cemoqt}.

\section*{Acknowledgments} 
ANWH is supported by 
Fellowship EP/M004333/1  from the 
Engineering \& Physical Sciences Research Council, UK, and he thanks the 
School of Mathematics and Statistics, University of New South Wales, for 
hosting him as a Visiting Professorial Fellow with funding from the Distinguished Researcher Visitor Scheme. 
He is also grateful to   
John Roberts and Wolfgang Schief for providing additional support during his 
stay in Sydney, where the idea behind this work originated, and to his colleagues at La Trobe 
for hospitality 
during his visit to  Melbourne in May 2019.
Part of this work was carried out at the Isaac Newton Institute, during the GCS programme, supported by EPSRC grant EP/R014604/1. GRWQ is grateful to the Simons Foundation for a grant supporting this work. 

\end{document}